\documentclass[11pt]{amsart}

\usepackage{amsfonts}
\usepackage{mathrsfs}
\usepackage{bbm}
\usepackage{amsmath}
\usepackage{amssymb}
\usepackage{amscd}
\usepackage{amsthm}

\makeatletter
\@namedef{subjclassname@2020}{\textup{2020} Mathematics Subject Classification}
\makeatother

\newcommand{\di}{\displaystyle}

\newtheorem{thm}{Theorem}[section]

\newtheorem{prop}[thm]{Proposition}
\newtheorem{cor}[thm]{Corollary}

\newcommand{\bc}{\mathbb{C}}

\newcommand{\bp}{\mathbb{P}}

\newcommand{\scc}{\mathscr{C}}
\newcommand{\sd}{\mathrm{od}}
\newcommand{\scp}{\mathscr{P}}
\newcommand{\st}{\mathscr{T}}
\newcommand{\sj}{\mathscr{J}}

\newcommand{\fl}{\mathfrak{L}}
\newcommand{\fp}{\mathfrak{P}}

\newcommand{\fs}{\mathfrak{S}}

\newcommand{\fsi}{\mathfrak{S}_{\infty}}
\newcommand{\fsz}{\mathfrak{S}_{0}}

\newcommand{\co}{\mathcal{O}}

\newcommand{\cisom}{C^{\infty}\left(S^1, M\right)}

\newcommand{\cisopo}{C^{\infty}\left(S^1, \mathbb{P}^1\right)}

\newcommand{\cP}{\Lambda}

\newcommand{\ltz}{\Lambda_{t_0}^n}

\newcommand{\lbn}{\mathfrak{L}_n}
\newcommand{\lbnk}{\mathfrak{L}_{n, k}}

\newcommand{\R}{\mathrm{R}(\mathfrak{H}_{\varphi})}

\newcommand{\jnmpo}{J^{n-1}\left(t_0, \mathbb{P}^1\right)}
\newcommand{\jkpo}{J^k(t_0, \mathbb{P}^1)}

\newcommand{\jnmc}{J^{n-1}\left(t_0, \mathbb{C}\right)}
\newcommand{\jnc}{J^{n}\left(t_0, \mathbb{C}\right)}
\newcommand{\jkc}{J^{k}\left(t_0, \mathbb{C}\right)}

\newcommand{\jmz}{J^{n-1}\left(t_0, \mathbb{P}^1 \setminus \{0\} \right)}

\newcommand{\jnm}{j^{n-1}_{t_0}}
\newcommand{\jk}{j^{k}_{t_0}}

\newcommand{\jkx}{j^{k}_{t_0} x}

\newcommand{\sjt}{\mathscr{J}_{t_0}}

\newcommand{\sjm}{\mathscr{J}^{n-1}_{t_0}}
\newcommand{\sjmx}{\mathscr{J}^{n-1}_{t_0} x}

\newcommand{\cisoc}{C^{\infty}\left(S^1, \mathbb{C} \right)}
\newcommand{\cisoz}{C^{\infty}\left(S^1, \mathbb{P}^1 \setminus \{0\} \right)}

\newcommand{\cisoca}{C^{\infty}\left(S^1, \mathbb{C}^{\ast} \right)}

\newcommand{\hocipo}{\di H^0\left(C^{\infty}(S^1, \mathbb{P}^1), \ltz \right)}

\newcommand{\hojnm}{\di H^0\left(J^{n-1}\left(t_0, \mathbb{P}^1\right), \lbn \right)}

\newcommand{\hojk}{H^0\left(J^{k}\left(t_0, \mathbb{\bp}^1\right), \lbnk \right)}

\newcommand{\pn}{\mathscr{P}_{n}}
\newcommand{\pnd}{\mathscr{P}_{n, d}}

\newcommand{\pnl}{\mathscr{P}_{n, 0}}
\newcommand{\pny}{\mathscr{P}_{n, 1}}
\newcommand{\pnt}{\mathscr{P}_{n, 2}}

\newcommand{\pnn}{\mathscr{P}_{n, n}}

\newcommand{\pndl}{\mathscr{P}_{n, d, l}}
\newcommand{\pntl}{\mathscr{P}_{n, 2, l}}

\newcommand{\rto}{D^{t_0} \rho}

\begin{document}

\title{Holomorphic sections of line bundles on the jet spaces of the Riemann sphere}

\author{Xiaokun Wang}

\address{School of Information Science and Engineering, Shandong Univ., Qingdao, 266237,  P. R. China}

\email{wangxiaokun@sdu.edu.cn}

\author{Ning Zhang}

\address{School of Mathematical Sciences, Ocean Univ. of China, Qingdao, 266100, P. R. China}

\email{nzhang@ouc.edu.cn}



\thanks{The authors are grateful to L. Lempert for his helpful comments on the manuscript.
This research was partially supported by the Scientific Research Foundation
of Ocean University of China grant 861701013110.}


\begin{abstract}
  Fix a point $t_0$ in the circle $S^1$. The space $\jkpo$ of $k$-jets at $t_0$ of $C^{\infty}$ maps from
  $S^1$ to the Riemann sphere $\bp^1$ is a $k+1$ dimensional complex algebraic manifold.
  We identify a class of holomorphic sections of line bundles on $\jkpo$.
\end{abstract}

\maketitle

\begin{center}
\it Dedicated to Professor L\'{a}szl\'{o} Lempert on his seventieth birthday
\end{center}

\pagestyle{myheadings} \markboth{\centerline{X. Wang and N. Zhang}}{\centerline{Holomorphic sections of line bundles on the jet spaces}}


\section{Introduction \label{intro}}

Let $M$ be a positive dimensional complex manifold without boundary,
$t_0 \in S^1$ and $J^k(t_0, M)$, where $k=0, 1, 2, \cdots$, the submanifold of the $k$-jet space $J^k(S^1, M)$
consisting of $k$-jets at $t_0$.
Then the loop space $\cisom$ of all $C^{\infty}$ maps $S^1 \to M$
is a complex Fr\'echet manifold (see \cite{l04, ls}),
$J^k(t_0, M)$ is a finite dimensional complex manifold,
the map $$j^k_{t_0}: \cisom \ni x \mapsto j^k_{t_0} x \in J^k(t_0, M)$$ is holomorphic
and the target map
$$\tau_k: J^k(t_0, M) \ni j_{t_0}^k x \mapsto  x(t_0) \in M$$
gives rise to a holomorphic fiber bundle. Note that $j^0_{t_0}$
is simply the evaluation at $t_0$
$$E_{t_0}: \cisom \ni x \mapsto x(t_0) \in M.$$
If $M$ is complex algebraic, then both $J^k(t_0, M)$
and $\tau_k$ are algebraic. For example, $J^k(t_0, \bp^1)$
is a complex algebraic manifold of dimension $k+1$.
It follows from \cite[Section 4]{z03} that a generic complex algebraic hypersurface
of $\jkpo$ is not of the form $\tau_k^{-1}(W)$, where $W \subset \bp^1$.
In this paper, we study holomorphic sections of line bundles over
$\jkpo$.

Recall the group $PGL(2, \bc)$ of holomorphic automorphisms of $\bp^1$. The action of $PGL(2, \bc)$
on $\bp^1$ induces holomorphic actions
\begin{eqnarray*}
 & PGL(2, \bc) \times \cisopo \ni (\gamma, x) \mapsto \gamma \circ x \in \cisopo \hspace{2mm} \text{and} & \\
 & PGL(2, \bc) \times \jkpo \ni (\gamma, \jkx) \mapsto \jk (\gamma \circ x) \in \jkpo &
\end{eqnarray*}
with which the maps $\di j^k_{t_0}$ and $\tau_k$ are $PGL(2, \bc)$-equivariant.
The group $PGL(2, \bc)$ also acts on the Picard group of holomorphic line bundles on $\cisopo$ (resp. on $\jkpo$)
by pullbacks. Let $\di \ltz$ be the pullback line bundle $\di E_{t_0}^{\ast} \co(n)$ over $\cisopo$.
It follows from (\cite{z03, lz, z10}) that the group of $PGL(2, \bc)$-invariant holomorphic line bundles on $\cisopo$
is an infinite dimensional Lie group;
any $PGL(2, \bc)$-invariant line bundle on $\cisopo$ with non-zero holomorphic sections is of the form
$$\cP=\cP_{t_1}^{n_1} \otimes \cdots \otimes
\cP_{t_r}^{n_r},$$
where $n_i \ge 0$ and $t_i \neq t_j$ for $i\neq
j$; the space $\di H^0(\cisopo,\cP)$ of holomorphic sections of $\cP$ is finite dimensional and if $\cP$ is trivial,
then $\dim H^0(\cisopo,\cP)=1$.
Let $\di \lbnk$ be the pullback line bundle $\di \tau_{k}^{\ast} \co(n)$ over $\jkpo$.
Note that $\di E_{t_0}=\tau_{k}  \circ \jk$. So
\begin{equation*} \label{pb}
\ltz=\left(\jk \right)^{\ast} \lbnk.
\end{equation*}
In particular, every holomorphic function on $\jkpo$ is constant.

It follows from \cite[Section 4]{z03} that
for any positive integer $n$, we have
\begin{equation} \label{pullback}
\di \hocipo=\left(\jk \right)^{\ast} H^0\left(J^{k}\left(t_0, \mathbb{\bp}^1\right), \lbnk \right), \hspace{2mm} k \ge n-1;
\end{equation}
and if $n \ge 2$, $1 \le k \le n-1$, then
\begin{equation} \label{nmt}
\left(j_{t_{0}}^{k-1}\right)^{\ast} H^{0}\left(J^{k-1}\left(t_{0}, \mathbb{P}^{1}\right), \fl_{n, k-1} \right)
\subsetneqq \left(j_{t_{0}}^{k}\right)^{\ast} H^{0}\left(J^{k}\left(t_{0}, \mathbb{P}^{1}\right), \lbnk \right).
\end{equation}
Fix $x_0 \in \cisopo$. Define
    $$ Z=Z(\infty, S^1, \bp^1, t_0, x_0)=\{x \in \cisopo: j^{\infty}_{t_0} x=j^{\infty}_{t_0} x_0 \}. $$
    This is a connected complex submanifold of $\cisopo$ and any holomorphic
function on it is constant, see (3.1) and Theorem 4.2 of \cite{l04}.
Let $\fl$ be a holomorphic line bundle on $\jkpo$. The pullback bundle $\left(j^k_{t_0} \right)^{\ast} \fl$
is trivial on $Z$. If there exists a non-zero holomorphic section $\sigma$ of the bundle $\left(j^k_{t_0} \right)^{\ast} \fl$, 
then for any $x \in \cisopo$, $\sigma(x)$ only depends on the jets of $x$ at $t_0$.
So any $PGL(2, \bc)$-invariant line bundle on $\jkpo$ with non-zero holomorphic sections is of the form $\lbnk$, where $n \ge 0$.
We write $\di \lbn$ for $\fl_{n, n-1}$ ($n \ge 1$).
It turns out that the spaces
$\di H^{0}\left(J^{k}\left(t_{0}, \mathbb{P}^{1}\right), \lbnk \right)$, where $k=0, 1, \cdots$,
are completely determined by  $\di \hojnm$  (see Section \ref{bg}).

Fix a holomorphic section $\di \fsi=\fs_{\infty, n}\in H^0\left(\bp^1, \co(n) \right)$ such that the only zero of $\fsi$ is $\infty \in \bp^1$.
The pullback section
$$\tau_{n-1}^{\ast} \fsi \in \hojnm$$
does not vanish on the dense open subset
$\di \jnmc \subset \jnmpo$. Let $\pn$ be the range of the linear operator
$$\hojnm \ni \varsigma \mapsto \left. \frac{\varsigma}{\tau_{n-1}^{\ast} \fsi}\right|_{\jnmc} \in \co\left(\jnmc \right).$$
It follows from \cite[Section 4]{z03} and (\ref{pullback}) that any $\rho \in \pn$ is a
polynomial of degree less than or equal to $n$ in  $\di x(t_0)$, $\di x^{(1)}(t_0)$, $\cdots$, $\di x^{(n-1)}(t_0)$
and the degree of $\rho$ as a polynomial of $\di x^{(n-1)}(t_0)$
is smaller than or equal to $1$.
If we consider elements of $\scp_{n-1}$ as polynomials on $\jnmc$ independent
of $\di x^{n-1}(t_0)$, then
\begin{equation*} \label{inclusion}
\scp_{n-1} \subset \pn.
\end{equation*}
Let $\scp_{n, d}$ be the subspace of $\pn$ consisting of homogeneous
polynomials of degree $d$, where $d=0, 1, \cdots, n$. Then
\begin{eqnarray*}
& \pn=\oplus_{d=0}^{n} \pnd, & \\
& \dim \pnd=\dim \scp_{n, n-d}, &
\end{eqnarray*}
$\di \dim \scp_{n, 0}=1$ and the set of polynomials
$\di x(t_0)$, $\di x^{(1)}(t_0)$, $\cdots$, $\di x^{(n-1)}(t_0)$
is a basis for $\pny$ (for all above see Section \ref{bg}).

Given a monomial $\di \rho_{\mu}=c_{\mu} x^{(r_1)}(t_0) \cdots x^{(r_m)}(t_0)$, where $c_{\mu} \not=0$, we call $r_1+\cdots+r_m$
the order of derivatives of $\rho_{\mu}$ and denote it by $\sd(\rho_{\mu})$.
In particular, the order of derivatives of a non-zero constant is $0$.
For a non-zero polynomial $\rho=\sum_{\mu} \rho_{\mu}$  in finitely many  derivatives $\di x^{(\nu)}(t_0)$, we define
$$\sd(\rho)=\max_{\mu} \left\{\sd \left(\rho_{\mu} \right) \right\}.$$
Let $\pndl \subset \pnd$ be the subspace spanned by polynomials $\rho \in \pnd$ such
that the order of derivatives of every monomial in $\rho$ is $l$.
Then $$\di \dim \pndl=\dim \scp_{n, n-d, l}$$ (see Section \ref{ci}).
Our main result is

\begin{thm} \label{main}
Let $n \ge 2$ be an integer. Then
the following set of polynomials is a basis for $\pnt$:
\begin{equation} \label{degreetwo}
x^{(r_1)}(t_0)x^{(r_2)}(t_0), \hspace{2mm} x^{(k_1)}(t_0)x^{(k_2)}(t_0)-\frac{k_1}{k_2+1}x^{(k_1-1)}(t_0)x^{(k_2+1)}(t_0),
\end{equation}
where $\di 0 \le r_2 \le r_1\le n-2-r_2$ and $\di 0 \le n-1-k_1 \le k_2 \le k_1-2$.
In particular,
$$\di \pnt=\oplus_{l=0}^{2n-4} \scp_{n, 2, l} \hspace{2mm} \text{and} \hspace{2mm} \dim \pnt=\binom{n}{2}.$$
\end{thm}

It is straightforward to verify that
\begin{equation*}
\dim \pntl=\left\{\begin{array}{ll}
   \di 1+\left[\frac{l}{2} \right], & 0 \le l \le n-2,\\
   \vspace{-3mm} & \\
   \di 1+\left[\frac{2n-4-l}{2} \right], & n-1 \le l \le 2n-4
\end{array} \right.
\end{equation*}
(where $\left[l/2 \right]$ is the integer part of $l/2$).
So $1 \le \dim \pntl \le \left[\frac{n}{2} \right]$.

Note that for $n=2, 3, \cdots$ and $d=0, 1, 2, n-2, n-1, n$, we have
\begin{itemize}
  \item[(i)] $\dim \pnd=\binom{n}{d}$;

  \item[(ii)] $\di \pnd=\oplus_{l=0}^{dn-d^2} \scp_{n, d, l}$, where $\dim \scp_{n, d, l} \not=0$ for $0 \le l \le dn-d^2$;

  \item[(iii)] $\dim \scp_{n, d, l}=\dim \scp_{n, d, dn-d^2-l}.$
\end{itemize}
It is natural to ask whether (i), (ii) and (iii) are still true for $n=6, 7, \cdots$ and $d=3$, $4$, $\cdots$, $n-3$.

This paper is organized as follows. In Section \ref{bg}, we recall some of the results in \cite{z03}.
A polynomial $\rho_1$ on $\jnmc$ is in $\pn$ if and only if there is a polynomial $\rho_2$ on $\jnmc$
such that $\rho_1$ and $\rho_2$ satisfy certain compatibility condition.
In Section \ref{ci}, we study these compatible pairs of polynomials.
It turns out that any polynomial in (\ref{degreetwo}) consisting of more than one monomials
is closely related to combinatorial identities.
In Section \ref{combi}, we prove a family of combinatorial identities (see Proposition \ref{cis}).
In the final Section \ref{final}, we prove Theorem \ref{main}. As a consequence of Theorem \ref{main},
we obtain more combinatorial results (see Corollary \ref{general}).


\section{Background \label{bg}}

 Let $\di \fsi \in H^0\left(\bp^1, \co(n) \right)$ be as in Section \ref{intro}.
 The pullback section $$\di E_{t_0}^{\ast} \fsi \in \hocipo$$ does not vanish on the dense open
subset $\cisoc \subset \cisopo$.
We write $\R$ for the range of the linear operator
$$\hocipo \ni \sigma \mapsto \left. \frac{\sigma}{E_{t_0}^{\ast} \fsi} \right|_{\cisoc} \in \co \left(\cisoc \right).$$
It follows from \cite[Section 4]{z03} (where we set $\varphi=E^n_{t_0}$ in \cite[(4.1)]{z03}) that every element $P$ of $\R$
is a polynomial in the linear functionals $\di x(t_0)$, $\di x^{(1)}(t_0)$, $\cdots$, $\di x^{(n-1)}(t_0)$
(see Proposition 4.2 and Theorem 4.7);
the degree of $P$ as a polynomial of $\di x^{(n-1)}(t_0)$
is smaller than or equal to $1$ (see Proposition 4.3); and $x(t_0), x^{(1)}(t_0), \cdots, x^{(n-1)}(t_0) \in \R$
(see the proof of Theorem 4.7).
The bundle $\ltz$ is $PGL(2, \bc)$-invariant. So we have (\ref{pullback}) (here we follow \cite[(31)]{z21}). In particular,
\begin{equation*}
  H^0\left(J^k(t_0, \bp^1), \fl_{1, k} \right)=\tau_k^{\ast} H^0\left(\bp^1, \co(1) \right), \hspace{2mm} k=1, 2, \cdots.
\end{equation*}
It is clear that the set of polynomials $x(t_0), x^{(1)}(t_0), \cdots, x^{(n-1)}(t_0)$
is a basis for $\pny$.
Let $\di \varsigma \in \hojk$, where $k=1, 2, \cdots, n-2$. Recall that $\ltz=\left(\jk \right)^{\ast} \lbnk$. So the function
\begin{equation*}
  \left. \frac{\varsigma}{\tau_k^{\ast} \fsi} \right|_{\jkc}
\end{equation*}
is a polynomial in $\di x(t_0), x^{(1)}(t_0), \cdots, x^{(k)}(t_0)$, which can be considered as an element of $\pn$
independent of the variables $x^{(k+1)}(t_0)$, $\cdots$, $x^{(n-1)}(t_0)$. Thus we have (\ref{nmt}).

The manifold $\jnmpo$ is covered by open subsets $\jnmc$ and $\jmz$.
For any $x \in \cisopo$ with $x(t_0) \in \bc$, we write $\sjmx$ for the vector
$$\left(x(t_0), x^{(1)}(t_0), \cdots, x^{(n-1)}(t_0) \right) \in \bc^{n}.$$
Then the maps
\begin{equation*}
  \jnmc \ni \jnm x \mapsto \sjmx \in \bc^n \hspace{2mm} \text{and}
\end{equation*}
\begin{equation*}
  \jmz \ni \jnm x \mapsto \sjm x^{-1} \in \bc^n
\end{equation*}
give rise to local charts on $\jnmpo$.
Let $\rho_1, \rho_2$ be holomorphic functions on $\jnmc$ and let $\di \bc^{\ast}=\bc \setminus \{0\}$.
We write $\rho_1 \sim_n \rho_2$ if
\begin{equation} \label{compatible}
\rho_1\left( \sjmx \right)=x(t_0)^n \rho_2\left( \sjm x^{-1} \right), \hspace{3mm} x \in \cisoca.
\end{equation}
The condition (\ref{compatible}) implies that $\rho_2 \sim_n \rho_1$.
Let $\fsz \in H^0\left(\bp^1, \co(n) \right)$ be the section with $\di \left. \left(\fsz/\fsi\right) \right|_{\bc}=\zeta^n$, where $\zeta \in \bc$.
The bundle $\lbn$ is $PGL(2, \bc)$-invariant.
Thus for any $\varsigma \in \hojnm$,
 there are polynomials $\rho_1=\rho_1(\varsigma)$ and $\rho_2=\rho_2(\varsigma)$ on $\jnmc \simeq \bc^n$ with $\rho_1 \sim_n \rho_2$ such that
\begin{eqnarray*}
  \rho_1\left( \sjmx \right) &=& \frac{\varsigma}{\tau_{n-1}^{\ast} \fsi}\left(\jnm x \right), \hspace{2mm} x \in \cisoc, \\
  \rho_2\left( \sjm x^{-1} \right) &=& \frac{\varsigma}{\tau_{n-1}^{\ast} \fsz}\left(\jnm x \right), \hspace{2mm} x \in \cisoz.
\end{eqnarray*}
It is clear that a pair of holomorphic functions $\rho_1, \rho_2$ on $\jnmc$ with  $\rho_1 \sim_n \rho_2$
determine a section $\varsigma \in \hojnm$. So $\rho_1, \rho_2$ must be polynomials; and we obtain the following

\begin{prop} \label{iff}
  A holomorphic function
$\rho_1$ on $\jnmc$ is in $\pn$ if and only if there is a holomorphic function $\rho_2$ on $\jnmc$ such that $\rho_1 \sim_n \rho_2$.
\end{prop}

Suppose $\di \rho_1 \sim_n \rho_2$. Let $\rho_i^d$ be the $d$-th order homogeneous
component of $\rho_i$, $i=1, 2$, $d=0, 1, \cdots, n$. It follows from \cite[Proposition 3.2]{z03}
and (\ref{pullback}) that $$\di \rho_1^d \sim_n \rho_2^{n-d}.$$
So $\di \pn=\oplus_{d=0}^{n} \pnd$ and
$\di \dim \pnd=\dim \scp_{n, n-d}.$
In particular, $\di \dim$ $\scp_{n, n-1}$ $=$ $n$ and $\dim \pnl=\dim \pnn=1$ ($1 \sim_n x(t_0)^n$).
If $\rho_1 \sim_{n-1} \rho_2$, then $\di \rho_1 \sim_n x(t_0) \rho_2$.
Hence $\scp_{n-1} \subset \pn$.

Let $\rho$ be a polynomial in $\di \sjmx$, where $x \in \cisoc$.
Letting $t_0$ vary in $S^1$, $\rho$ induces a $C^{\infty}$ function
$$\chi: \cisoc \times S^1 \ni (x, t_0) \mapsto \rho\left(\sjmx \right) \in \bc.$$
We define a linear operator $D^{t_0}$ from the space of polynomials on $\jnmc$ to
the space of polynomials on $\jnc$ by
\begin{equation*}
  \rto\left(\sjt^n x \right)=\frac{\partial}{\partial t_0} \chi(x, t_0)
\end{equation*}
(for example, if $\rho=x^{(r)}(t_0)$, then $\di \rto=x^{(r+1)}(t_0)$).
It follows from \cite[Proposition 4.6]{z03} that $D^{t_0}$ maps $\pnd$ into $\scp_{n+1, d}$.


\section{Compatible polynomials $\rho_1 \sim_n \rho_2$ \label{ci}}

Let $m, k$ be positive integers, $\fp_{m, k}$ the space of $m$-tuples of nonnegative integers
$\di p=(p_1, \cdots, p_m)$ with $$\di p_1+2p_2+\cdots +m p_m=k$$ and $|p|=p_1+\cdots+p_m$.
We write $\fp_k$ for the space $\fp_{k, k}$.
Given $\di x \in \cisopo$ with $x(t_0) \not=0$, it follows from the Fa\`{a} di Bruno's formula that
\begin{equation} \label{chain}
\left(\frac{1}{x} \right)^{(k)}(t_0)  =  \sum_{p \in \fp_k}  c_p \frac{x^{(1)}(t_0)^{p_{1}}
\cdots x^{(k)}(t_0)^{p_{k}}}{x(t_0)^{|p|+1}},
\end{equation}
where
\begin{equation*}
c_p=k ! \frac{(-1)^{|p|}}{(1!)^{p_1}  \cdots (k!)^{p_k}} \binom{|p|}{p_1, \cdots, p_k} \hspace{2mm}
\text{and} \hspace{2mm} \binom{|p|}{p_1, \cdots, p_k}=\frac{|p|!}{p_1! \cdots p_k!}.
\end{equation*}
On the right hand side of (\ref{chain}), we have $1 \le |p| \le k$;
the only term with $|p|=1$ resp. $|p|=k$ is
$$\di \frac{-x^{(k)}(t_0)}{x(t_0)^2}
\hspace{3mm} \text{resp.} \hspace{3mm}  \frac{(-1)^k k! x^{(1)}(t_0)^k}{x(t_0)^{k+1}}$$
and the order of derivatives of the numerator of every term is $k$.

Let $k_1 \ge k_2$ be positive integers and $\di \mu=(\mu_1, \cdots, \mu_{k_1}) \in \fp_{k_1, k_1+k_2}$. If
there exist
$\di p=(p_1, \cdots, p_{k_1}) \in \fp_{k_1}$ and
$\di q=(q_1, \cdots, q_{k_2}) \in \fp_{k_2}$ such that
$$\mu=p+q \triangleq (p_1+q_1, \cdots, p_{k_2}+q_{k_2}, p_{k_2+1}, \cdots, p_{k_1}),$$
then we define
\begin{equation*}
  \scc_{\mu, k_1, k_2}= \sum_{p \in \fp_{k_1}, q \in \fp_{k_2}, p+q=\mu}
  \binom{|p|}{p_1, \cdots, p_{k_1}} \binom{|q|}{q_1, \cdots, q_{k_2}}.
\end{equation*}
Otherwise we set $\di  \scc_{\mu, k_1, k_2}=0$.
By (\ref{chain}), we have
\begin{equation} \label{twoterms}
  \left(\frac{1}{x} \right)^{(k_1)}(t_0) \left(\frac{1}{x} \right)^{(k_2)}(t_0) =\sum_{\mu \in \fp_{k_1, k_1+k_2}} c_{\mu} \frac{x^{(1)}(t_0)^{\mu_{1}}
  \cdots x^{(k_1)}(t_0)^{\mu_{k_1}}}{x(t_0)^{|\mu|+2}},
\end{equation}
where
$$c_{\mu}=(k_1)!(k_2)!   \frac{(-1)^{|\mu|}}{(1!)^{\mu_1}  \cdots (k_1!)^{\mu_{k_1}}} \scc_{\mu, k_1, k_2}.$$

Similarly, for any monomial $\rho_{\mu} \not=0$ in $x(t_0)$, $x^{(1)}(t_0)$, $\cdots$, $x^{(m)}(t_0)$,
the function $\di \rho_{\mu} \left( \sj_{t_0}^m x^{-1} \right)$ ($x(t_0) \not=0$) can be expressed
as a finite linear combination of pairwise distinct rational functions in $x(t_0)$, $x^{(1)}(t_0)$, $\cdots$, $x^{(m)}(t_0)$, where
the numerator of each rational function is
a monic monomial in $x^{(1)}(t_0)$, $\cdots$, $x^{(m)}(t_0)$ whose order of derivatives is  $\sd(\rho_{\mu})$ and
the denominator is a power of $x(t_0)$.
We call this linear combination the inversion expansion of $\di \rho_{\mu}$.
The inversion expansion of a polynomial $\rho =\sum_{\mu} \rho_{\mu}$
is obtained by taking the sum of the inversion expansions of the monomials $\rho_{\mu}$.
If $\rho_1 \sim_n \rho_2$, then $\sd(\rho_1)=\sd(\rho_2)$. So $\di \dim \pndl=\dim \scp_{n, n-d, l}$.
By Proposition \ref{iff}, we have the following

\begin{prop} \label{monoj}
    Let $\rho$ be a polynomial in $x(t_0)$, $x^{(1)}(t_0)$, $\cdots$, $x^{(m)}(t_0)$.
    Then $\rho \in \scp_{n} \setminus \scp_{n-1}$ if and only if
    in the inversion expansion of $\rho$, the highest power of $x(t_0)$ in the denominators is $x(t_0)^{n}$.
    In particular, the monomial
    \begin{equation*}
  x^{(k_1)}(t_0) \cdots x^{(k_d)}(t_0) \in \scp_{k_1+\cdots+k_d+d} \setminus \scp_{k_1+\cdots+k_d+d-1},
    \end{equation*}
    where $d=1, 2, \cdots$ and $k_1, \cdots, k_d$ are nonnegative integers.
\end{prop}

\begin{prop} \label{sameTOD}
  Let $\rho =\sum_{\mu} \rho_{\mu} \in \pnd$, where $\rho_{\mu}$ are monomials, and let $\rho_{\mu_0}$ be a fixed monomial in $\rho$. Then
  $$\di \sum_{\sd(\rho_{\mu})=\sd(\rho_{\mu_0})} \rho_{\mu} \in \pnd,$$
  where the summation is taken over all monomials $\rho_{\mu}$ in $\rho$ with $\sd(\rho_{\mu})=\sd(\rho_{\mu_0})$.
\end{prop}
\begin{proof}
Recall that the order of derivatives of the numerator of every term in the inversion
expansion of $\rho_{\mu}$ is $\sd(\rho_{\mu})$. Note that
\begin{eqnarray*}
  \di x(t_0)^n \rho\left( \sjm x^{-1} \right) &=& x(t_0)^n\sum_{\sd(\rho_{\mu}) = \sd(\rho_{\mu_0})} \rho_{\mu}\left( \sjm x^{-1} \right) \\
  & + & x(t_0)^n \sum_{\sd(\rho_{\mu}) \not= \sd(\rho_{\mu_0})} \rho_{\mu}\left( \sjm x^{-1} \right)
\end{eqnarray*}
is a polynomial in $\di \sjmx$, where $x \in \cisoca$. So the first part on the right hand side of the above equation
is a polynomial in $\di \sjmx$, from which the conclusion of the proposition follows.
\end{proof}

By Proposition \ref{sameTOD}, $\di \pnd$ is the direct sum of finitely many subspaces of the form $\di \pndl$.


\section{Combinatorial identities \label{combi}}

Any $r=(r_1, \cdots, r_m) \in \fp_{m, k}$ induces a multiset of positive integers $$\di M_r=\{r_1 \cdot 1, r_2 \cdot 2, \cdots, r_{m} \cdot m \}$$
  (where $r_1, \cdots, r_m$ are the repetition numbers). Recall that the number of permutations of $M_r$ is $\binom{|r|}{r_1, \cdots, r_m}$.
  Let
  $$\di \omega=(a_1, \cdots, a_{|r|})$$
  be a permutation of $M_{r}$. The $i$-th partial sum $S_i=S_i(\omega)$ of $\omega$ is the sum of the first $i$ terms of $\omega$.
  The sequence of partial sums $S_1, \cdots, S_{|r|}$ of $\omega$ is strictly increasing and $\di S_{|r|}=k$.
  Note that for any $\mu \in \fp_{n-1, 2n-4}$ with $|\mu| \ge n-1$, where $n \ge 4$,
  we must have $\mu_{n-1}=0$. Thus we may also consider $\mu$ as an element of $\di \fp_{n-2, 2n-4}$.
  Recall the number $\di  \scc_{\mu, k_1, k_2}$ defined in Section \ref{ci}.

\begin{prop} \label{cis}
  Let $n \ge 4$ be an integer. Then for any $\mu \in \fp_{n-1, 2n-4}$ with $|\mu| \ge n-1$, we have
  $$\di  \scc_{\mu, n-1, n-3}=\di  \scc_{\mu, n-2, n-2}.$$
\end{prop}
\begin{proof}
Let $A$ resp. $A'$ be the set of permutations $\omega$ of the multiset $M_{\mu}$ such that
the sequence of partial sums of $\omega$ contains
$n-1$ resp. $n-2$.
  Any element of $A$ is of the form
  $$\omega=(a_1, \cdots, a_{|p|}, b_1, \cdots, b_{|q|}),$$
  where $\omega_1=(a_1, \cdots, a_{|p|})$ is a permutation of the multiset $M_{p}$ for some $p \in \fp_{n-2, n-1}$,
  $\omega_2=(b_1, \cdots, b_{|q|})$ is a permutation of the multiset $M_{q}$ for some $q \in \fp_{n-3}$ and $\mu=p+q$.
  Let $S_1, \cdots, S_{|p|}$ be the sequence
  of partial sums of $\omega_1$ and $\tilde{S}_1, \cdots, \tilde{S}_{|q|}$ the sequence
  of partial sums of $\omega_2$.
  Similarly, any element of $A'$ is of the form
  $$\omega'=(a_1, \cdots, a_{|\alpha|}, b_1, \cdots, b_{|\beta|}),$$
  where $\omega'_1=(a_1, \cdots, a_{|\alpha|})$ is a permutation of the multiset $M_{\alpha}$ for some $\alpha \in \fp_{n-2}$,
  $\omega'_2=(b_1, \cdots, b_{|\beta|})$ is a permutation of the multiset $M_{\beta}$ for some $\beta \in \fp_{n-2}$
  and $\mu=\alpha+\beta$.
  Let $S'_1, \cdots, S'_{|\alpha|}$ be the sequence
  of partial sums of $\omega'_1$ and $\tilde{S'}_1, \cdots, \tilde{S'}_{|\beta|}$ the sequence
  of partial sums of $\omega'_2$.

  Note that the number of elements in $A$ resp. $A'$ is $\di \scc_{\mu, n-1, n-3}$ resp.
 $\di  \scc_{\mu, n-2, n-2}$.
  Next we construct a bijection $\st: A \to A'$.
  If $a_1=1$, then we define $$\st(\omega)=(a_2, \cdots, a_{|p|}, 1, b_1, \cdots, b_{|q|}) \in A'.$$ If $a_1 >1$, then the set
  $$\{S_1-1, \cdots, S_{|p|}-1, \tilde{S}_1, \cdots, \tilde{S}_{|q|} \}$$
  contains $|\mu| \ge n-1$ integers between $1$ and $n-2$. Let $\di S_{i_1}$ be the smallest partial sum of $\omega_1$ such that
  there is a partial sum $\di \tilde{S}_{i_2}$ of $\omega_2$ with $\di S_{i_1}-1=\tilde{S}_{i_2}$
  (the existence of $\di S_{i_1}$ follows from the pigeonhole principle).
  Now we define
  $$\di \st(\omega)=(b_1, \cdots, b_{i_2}, a_{i_1+1}, \cdots, a_{|p|}, a_1, \cdots, a_{i_1}, b_{i_2+1}, \cdots, b_{|q|}) \in A'.$$
  Note that $\di \tilde{S}_{i_2}$ is the smallest partial sum of $\omega_2$ with $\di S_{i_3}-1=\tilde{S}_{i_2}$
  for some partial sum $\di S_{i_3}$ of $\omega_1$. It is clear that $\st$ is injective. Let $\omega' \in A'$.
  If $b_1=1$, then $\omega' \in \st(A)$. If $b_1>1$, applying the pigeonhole principle, one can find
  the smallest partial sum $\di \tilde{S'}_{i_4}$ of $\omega'_2$  with $\di \tilde{S'}_{i_4}-1=S'_{i_5}$
  for some partial sum $\di S'_{i_5}$ of $\omega'_1$. Thus we also have $\omega' \in \st(A)$.
\end{proof}


\section{The space $\pnt$ \label{final}}

\begin{prop} \label{maxTOD}
Let $n \ge 3$ be an integer. Then
  $$\rho_{n, 2, 2n-4}=\di x^{(n-1)}(t_0)x^{(n-3)}(t_0)-\frac{n-1}{n-2} x^{(n-2)}(t_0)^2 \in \pnt \setminus \scp_{n-1, 2}.$$
\end{prop}
\begin{proof} As $\di \rho_{n, 2, 2n-4}$ contains $\di x^{(n-1)}(t_0)$, we have $\di \rho_{n, 2, 2n-4} \not\in \scp_{n-1, 2}$.
A direct computation yields
$$\di -x^{(2)}(t_0) \sim_3 x^{(2)}(t_0)x(t_0)-2x^{(1)}(t_0)^2.$$
Thus $\di \rho_{3, 2, 2} \in \scp_{3, 2}$.
When $n \ge 4$, it follows from (\ref{twoterms}) and Proposition \ref{cis} that
in the inversion expansion of $\rho_{n, 2, 2n-4}$, the highest power of $x(t_0)$ in the denominators is $x(t_0)^{n}$.
So $\di \rho_{n, 2, 2n-4} \in \pnt$.
  \end{proof}

\vspace{2mm}

\noindent {\it Proof of Theorem \ref{main}.}
Note that $\scp_{2, 2}$ is spanned by $x(t_0)^2$.
Suppose that the conclusion of the theorem holds for $\scp_{n-1, 2}$, $n=3, 4, \cdots$. Next we
show that it is also true for $\pnt$.

First we verify that each polynomial in (\ref{degreetwo}) is an element of $\pnt$.
By Proposition \ref{monoj}, $\di x^{(r_1)}(t_0)x^{(r_2)}(t_0) \in \pnt$
if and only if $r_1+r_2 \le n-2$. Since $\scp_{n-1, 2} \subset \pnt$, we only need to prove that
\begin{equation} \label{nkt}
x^{(n-1)}(t_0)x^{(k_2)}(t_0)-\frac{n-1}{k_2+1}x^{(n-2)}(t_0)x^{(k_2+1)}(t_0) \in \pnt,
\end{equation}
where $\di 0 \le k_2 \le n-3$.
The case when $k_2=n-3$ follows from Proposition \ref{maxTOD}.
Applying the linear operator $D^{t_0}$ (see Section \ref{bg}) to
the polynomial
$$\di x^{(n-2)}(t_0)x^{(n-4)}(t_0)-\frac{n-2}{n-3} x^{(n-3)}(t_0)^2 \in \scp_{n-1, 2},$$
we obtain (\ref{nkt}) for $n \ge 4$ and $k_2=n-4$. Similarly, applying $D^{t_0}$ to
the polynomial
$$x^{(n-2)}(t_0)x^{(k_2)}(t_0)-\frac{n-2}{k_2+1}x^{(n-3)}(t_0)x^{(k_2+1)}(t_0) \in \scp_{n-1, 2},$$
where $n \ge 5$ and $\di 0 \le k_2 \le n-5$, we obtain a polynomial
\begin{eqnarray*}
  && \left( x^{(n-1)}(t_0)x^{(k_2)}(t_0)-\frac{n-1}{k_2+1}x^{(n-2)}(t_0)x^{(k_2+1)}(t_0) \right) \\
  &+& \frac{k_2+2}{k_2+1} \left( x^{(n-2)}(t_0)x^{(k_2+1)}(t_0)-\frac{n-2}{k_2+2}x^{(n-3)}(t_0)x^{(k_2+2)}(t_0) \right)
\end{eqnarray*}
in $\pnt$. The second part of the above polynomial is in $\di \scp_{n-1, 2} \subset \pnt$. So the first part is
also in $\pnt$.

Let $\rho \in \pnt$. Recall that any monomial in $\rho$ is of the form
\begin{equation*}
\di cx^{(i_1)}(t_0)x^{(i_2)}(t_0),
\end{equation*}
where $0 \le i_2 \le i_1 \le n-1$,
$i_1+i_2 \le 2n-3$ and $c$ is a constant. If $i_1+i_2=2n-3$, then we must have $i_1=n-1$ and $i_2=n-2$.
By Proposition \ref{sameTOD}, we have $\di x^{(n-1)}(t_0)x^{(n-2)}(t_0) \in \pnt$,
which is a contradiction to Proposition \ref{monoj}. Thus $i_1+i_2 \le 2n-4$
and $\di \pnt=\oplus_{l=0}^{2n-4} \scp_{n, 2, l}$.

Note that any polynomial $\rho$
in (\ref{degreetwo}) is an element of $\di \scp_{n, 2, \sd(\rho)}$.
Let $\rho_1$, $\cdots$, $\rho_m$ be pairwise distinct polynomials in (\ref{degreetwo}) with $\sd(\rho_1)=\cdots=\sd(\rho_m)$
and let $x^{(\nu_i)}(t_0)$ be the highest derivative that $\rho_i$ depends on. Then
$\nu_i \not=\nu_k$ if $i \not= k$. It is clear that $\rho_1$, $\cdots$, $\rho_m$ are linearly independent.

The space spanned by all polynomials in (\ref{degreetwo})
has the decomposition
$$\oplus_{l=0}^{2n-4} V_l, \hspace{3mm} \text{where} \hspace{2mm} V_l \subset \scp_{n, 2, l}.$$
Let $m_l$ be the number of monic monomials of degree $2$ on $\jnmc$ whose order of derivatives is $l$.
When $\di 0 \le l \le n-2$, we have $\dim V_l=m_l$.
So $\di V_l=\scp_{n, 2, l}$. If $n-1 \le l \le 2n-4$, then a single monomial cannot be an element of $\scp_{n, 2, l}$
(see Proposition \ref{monoj}).
Thus $\dim \scp_{n, 2, l}< m_l$.
As $\dim V_l=m_l-1$, we have $\di V_l=\scp_{n, 2, l}$.
\qed

\vspace{2mm}

  By considering the inversion expansion of the polynomial
\begin{equation*}
\di  x^{(n-1)}(t_0)x^{(k_2)}(t_0)-\frac{n-1}{k_2+1}x^{(n-2)}(t_0)x^{(k_2+1)}(t_0) \in \pnt \setminus \scp_{n-1, 2}
\end{equation*}
(see (\ref{twoterms})), we obtain the following generalization of Proposition \ref{cis}.

\begin{cor} \label{general}
   Let $n \ge 4$ and $1 \le k_2 \le n-3$ be integers. Then for any
   $\mu \in \fp_{n-1, n-1+k_2}$ with $|\mu| \ge n-1$, we have
   $$\scc_{\mu, n-1, k_2}=\scc_{\mu, n-2, k_2+1}.$$
\end{cor}

For any $\rho \in \pnt \setminus \scp_{n-1, 2}$,
there exists a non-zero term in the inversion expansion of $\rho$ whose denominator is $\di x(t_0)^n$ (see Proposition \ref{monoj}).
Thus the conclusion of Corollary \ref{general} fails for $|\mu|= n-2$.


\end{document}